\documentclass[12pt]{article}
\usepackage[expansion=false]{microtype}
\usepackage[reqno]{amsmath}
\usepackage{fixmath}
\usepackage{amssymb, amsthm, enumitem}
\usepackage{mathtools}
\usepackage{graphicx}

\usepackage{ifthen,amssymb,amsfonts,amsmath,amsthm,graphicx, color, tikz, pgf}

\usepackage{hyperref}

\newtheoremstyle{plainsl}%
    {\topsep}
    {\topsep}
    {\slshape} 
    {}
    {\normalfont\bfseries}
    {.}
    { }
    {}

\swapnumbers

{\theoremstyle{plainsl}
\newtheorem{thm}{Theorem}[section]
\newtheorem{lem}[thm]{Lemma}

\newtheorem{prop}[thm]{Proposition}}
{\theoremstyle{remark}
\newtheorem{ex}[thm]{Example}
}

\renewcommand\proof{\noindent\textsl{Proof. }}
\newcommand\sqr[2]{{\vbox{\hrule height.#2pt
   \hbox{\vrule width.#2pt height#1pt \kern#1pt
        \vrule width.#2pt}\hrule height.#2pt}}}
\renewcommand\qed{%
    \ifmmode\eqno\sqr53
    \else\nolinebreak\ \hfill\sqr53\medbreak\fi}

\numberwithin{equation}{section}

\newcommand\pmat[1]{\begin{pmatrix} #1 \end{pmatrix}}

\newcommand{\abs}[1]{\left | #1 \right |}
\newcommand{\lp}{\! \left (}
\newcommand{\rp}{\right )} 
   \newcommand{\lsb}{\left \{ }
\newcommand{\rsb}{\right \} }

\newcommand\ip[2]{\langle#1,#2\rangle}

\title{A Spectral Moore Bound for Bipartite Semiregular Graphs}

 \author{Sabrina Lato}

\begin{document}
\maketitle

\begin{abstract}
  Let \( b \lp k, \ell, \theta \rp \) be the maximum number of vertices of valency \( k \) in a \( \lp k, \ell \rp \)-semiregular bipartite graph with second-largest eigenvalue \( \theta \). We obtain an upper bound for \( b \lp k, \ell, \theta \rp \) for \( 0 < \theta < \sqrt{k-1}+\sqrt{\ell-1} \). This bound is tight when there exists a distance-biregular graph with particular parameters, and we develop the necessary properties of distance-biregular graphs to prove this.\end{abstract}

\section{Introduction}

For a graph \( G, \) let \( \lambda_1 \lp G \rp \geq \lambda_2 \lp G \rp \cdots \geq \lambda_{\min} \lp G \rp \) denote the eigenvalues of the adjacency matrix of \( G \) in nonincreasing order. 

The second-largest eigenvalue, \( \lambda_2 \lp G \rp \), plays a key role in the study of expanders, as discussed in papers such as those by Alon~\cite{alonExpanders}, Brouwer and Haemers~\cite{spectraGraphs}, and Hoory, Linial, and Widgerson~\cite{expanderApplications}. A result of Alon and Boppana~\cite{alonExpanders} and Serre~\cite{serre} implies that for any \( k \) and \( \theta < 2 \sqrt{k-1} \), there are only finitely many \( k \)-regular graphs \( G \) with \( \lambda_2 \lp G \rp \leq \theta \). It is natural to wonder about extremal examples, where the second eigenvalue is minimal for graphs of given valency and number of vertices, or where the number of vertices is maximal given the valency and second eigenvalue.

Nozaki~\cite{Nozaki2015} studied when, for a fixed valency and number of vertices, a graph has minimal second-largest eigenvalue. In order to better understand these graphs, he developed a linear programming bound on the number of vertices. This bound was subsequently used by Cioab\u{a}, Koolen, Nozaki, and Vermette~\cite{cioaba1} to bound the maximum number of vertices for a graph with given valency and second eigenvalue not exceeding some given \( \theta \). Equality holds precisely when there is a distance-regular graph of valency \( k \), second eigenvalue \( \theta \), and girth at least twice the diameter. For graphs with diameter at most six, there are infinitely many examples, but if the diameter is strictly greater than six, a result of Damerell and Georgiacodis~\cite{damerell} shows that no such graphs exist.

In a subsequent paper, Cioab\u{a}, Koolen, and Nozaki~\cite{bipartite} were able to exploit the structure of bipartite graphs to come up with a better bound on the maximum number of vertices of a regular bipartite graph with given second eigenvalue. Their bound is tight when there is a bipartite distance-regular graph with diameter \( d, \) \( d-1 \) distinct eigenvalues, and girth at least \( 2d-2. \) In the same paper, they show that such graphs must have diameter \( d < 15. \)

A third paper, by Cioab\u{a}, Koolen, Mimura, Nozaki, and Okuda~\cite{hypergraph}, found an upper bound for the maximum number of vertices of an \( r \)-regular, \( u \)-uniform hypergraph with given second eigenvalue. Since the eigenvalues of a hypergraph are closely related to the bipartite incidence graph, where the bipartition is defined by vertices and hyperedges, this result also gives an upper bound on the number of vertices of a semiregular bipartite graph with given valencies and second eigenvalue.

In this paper, we also consider the maximum number of vertices that a semiregular bipartite graph of given valencies and second eigenvalue can have. However, by generalizing the techniques used to prove the bound for bipartite graphs in~\cite{bipartite} to apply to semiregular bipartite graphs, we can obtain better bounds. Let \( B \lp k, \ell, 2t, c \rp \) be the \( 2t \times 2t \) tridiagonal matrix with lower diagonal \( \lp 1, \ldots, 1, c, \ell \rp \), zero along the main diagonal, and row sum alternating between \( k \) and \( \ell. \) That is,
\[ B \lp k, \ell, 2t, c \rp = \pmat{0 & k \\ 1 & 0 & \ell-1 \\ & 1 & 0 & k-1 \\ & & \ddots & \ddots & \ddots \\ & & & 1 & 0 & \ell-1 \\ & & & & c & 0 & k-c \\ & & & & & \ell & 0}. \]

If \( \theta \) is the second-largest eigenvalue of \( B \lp k, \ell, 2t, c \rp \), the maximum number \( b \lp k, \ell, \theta \rp \) of vertices of valency \( k \) in a \( \lp k, \ell \rp \)-semiregular graph with second eigenvalue at most \( \theta \) is
\[ 1 +  k \sum_{i=1}^{t-2} \lp \ell-1 \rp^i \lp k - 1 \rp^{i-1} + \frac{k \lp \ell -1 \rp^{t-1} \lp k-1 \rp^{t-2}}{c}. \]
We can similarly define \( B \lp k, \ell, 2t+1, c \rp \) to be the \( \lp 2t+1 \rp \times \lp 2t+1 \rp \) tridiagonal matrix with lower diagonal \( \lp 1, \ldots, c, \ell \rp \), and if \( \theta \) is the second-largest eigenvalue, then the maximum number of vertices of valency \( k \) in a \( \lp k, \ell \rp \)-semiregular graph with second eigenvalue at most \( \theta \) is
\[ \ell \sum_{i=0}^{t-2} \lp \ell-1 \rp^i \lp k - 1 \rp^{i} + \frac{\ell \lp \ell -1 \rp^{t-1} \lp k-1 \rp^{t-1}}{c}.\]

As discussed by Cioab\u{a}~\cite{spectralMooreOverview}, the bound for regular graphs can be viewed as a generalization of the Moore bound. In a similar way, the bounds for semiregular bipartite graphs based on the second eigenvalue can be viewed as a generalization of a semiregular Moore bound. Such bounds were studied by F\'abrega, Fiol, and Yebra~\cite{semiregularMoore}. More recently, Araujo-Pardo, Dalf\'o, Fiol, and L\'opez~\cite{araujopardo} improved these bounds, and, for certain choices of \( k \) and \( \ell \), constructed families of infinite \( \lp k, \ell \rp \)-semiregular graphs meeting these bounds. Feng and Li~\cite{fengLi} proved that any infinite family of \( \lp k, \ell \rp \)-semiregular graphs must have second eigenvalue at least \( \sqrt{k-1}+\sqrt{\ell-1}. \) Here, we consider graphs with second eigenvalue strictly smaller than \( \sqrt{k-1} + \sqrt{\ell-1} \) to improve the bound on the number of vertices.

Viewing hypergraphs as bipartite graphs, the bound when \( t \) is even is of the same form as the bound found in Cioab\u{a}, Koolen, Mimura, Nozaki, and Okuda~\cite{hypergraph}, although the values for \( c \) may not agree. The bound when \( t \) is odd is of a different form, and we will give an infinite family of examples of when our bound for semiregular bipartite graphs is tight, although the bound for hypergraphs is not.

Much in the way that the graphs must be distance-regular for the bounds in~\cite{cioaba1} or~\cite{bipartite} to be tight, we find that when the bounds for semiregular bipartite graphs are tight, the graphs are distance-biregular. This natural relaxation of the notion of distance-regular graphs was studied in papers by Delorme~\cite{delorme}, Fiol~\cite{pseudoGlobal}, and Godsil and Shawe-Taylor~\cite{distanceRegularised}. We will give a brief overview of these graphs, and develop further theory to show that a semiregular bipartite graph with diameter \( d \) and girth \( g \geq 2d-2 \) is distance-biregular. This generalizes the result of Abiad, van Dam, and Fiol~\cite{abiad} for bipartite regular graphs.

\section{Preliminaries}
Let \( G = \lp B, C \rp \) be a semiregular bipartite graph where the sets of the partition are \( B \) and \( C \). The \textit{biadjacency matrix} \( N \) from \( B \) to \( C \) is the \( \abs{B} \times \abs{C} \) matrix with rows indexed by vertices in \( B \) and columns indexed by vertices in \( C \), and the \( \lp b, c \rp \)-th entry is equal to 1 if \( b \) is adjacent to \( c \) and 0 otherwise. Note that the transpose of the biadjacency matrix from \( B \) to \( C \) is the biadjacency matrix from \( C \) to \( B \).

We can similarly define the \( i \)-th \textit{distance biadjacency matrix} \( N_i \) as the matrix with rows indexed by vertices in \( B \) and columns indexed by either the vertices of \( B \) or the vertices of \( C, \) depending on whether \( i \) is, respectively, even or odd, and with the \( \lp b, x \rp \)-th entry equal to 1 if \( b \) is at distance \( i \) from \( x \), and 0 otherwise.

Further, we may write the adjacency matrix \( A \) of \( G \) as a block matrix of the form
\[ A = \pmat{0 & N \\ N^T & 0}, \]
where \( N \) is the biadjacency matrix from \( B \) to \( C \).

Let \( M \) be a normal matrix. As described in Section 8.12 of Godsil and Royle~\cite{yellow}, \( M \) admits a spectral decomposition. This means that for every eigenvalue \( \theta, \) there is a corresponding spectral idempotent \( E_{\theta} \) such that
\[ M = \sum_{\theta} \theta E_{\theta}. \]
The spectral idempotents are in fact pairwise orthogonal, so for any polynomial \( p \), we have
\[ p \lp M \rp = \sum_{\theta} p \lp \theta \rp E_{\theta}. \]
We will be making use of the spectral decomposition for \( A \) and \( NN^T \).

A \textit{backtracking walk} is a walk that contains the subwalk \( aba \) for two vertices \( a \) and \( b \). A \textit{nonbacktracking walk} is a walk that is not backtracking.

An inner product between two polynomials can be defined with respect to some measure \( \mu \) by
\[ \ip{p}{q} := \int p \lp x \rp q \lp x \rp d \mu. \]

The sequence \( \lp p_i \rp_{i \geq 0} \) is a sequence of \textit{orthogonal polynomials} with respect to \( \mu \) if the polynomials form an orthogonal basis with regard to this inner product and \( p_i \) has degree \( i \) for \( i \geq 0. \) A theorem, usually attributed to Favard~\cite{favard}, states that a family of monic polynomials is orthogonal if and only if they satisfy a certain three-term recurrence. We will be making repeated use of the discrete version of this theorem, which can be found in texts such as Nikiforov, Suslov, and Uvarov~\cite{classicalOrthogonal}, to take advantage of some of the properties of orthogonal polynomials.

\begin{thm}[Nikiforov, Suslov, and Uvarov~\cite{classicalOrthogonal}]\label{favard}
  A sequence \( \lp p_i \lp x \rp \rp_{i \geq 0} \) of monic polynomials is orthogonal if and only if for all \( i \geq 1, \) there exist real numbers \( a_i, b_i \) such that
  \[ p_{i+1} \lp x \rp = \lp x -a_i \rp p_{i} \lp x \rp - b_{i-1} p_{i-1} \lp x \rp. \]
\end{thm}

We also wish to use a few properties of orthogonal polynomials. The next result is standard, and follows from the discussion on page 140 of Godsil~\cite{yellow}.

\begin{lem}\label{interlace} Let \( p_0, \ldots, p_n \) be a sequence of monic orthogonal polynomials and let \( \alpha \in \mathbb{R} \). Then the zeros of \( p_{n-1} \) strictly interlace the zeros of \( p_n + \alpha p_{n-1}. \)
\end{lem}

This and other standard properties of orthogonal polynomials allow us to derive the following result, which is proved, in a more general case, by Cohn and Kumar~\cite{cohn}.

\begin{prop}[Cohn and Kumar~\cite{cohn}]\label{posCoeff}Let \( p_0, \ldots, p_n \) be a sequence of monic orthogonal polynomials, let \( \alpha \in \mathbb{R} \), and let \( r_1 \) be the largest root of \( p_n \lp x \rp + \alpha p_{n-1} \lp x \rp. \) Then
  \[ \frac{p_n \lp x \rp + \alpha p_{n-1} \lp x \rp}{x-r_1} \]
  has positive coefficients in terms of \( p_0, \ldots, p_n. \)
\end{prop}

\section{Biregular trees}

We can define the \( \lp k, \ell \rp \)-biregular tree as the infinite tree where every vertex is either adjacent to \( k \) vertices of valency \( \ell \) or \( \ell \) vertices of valency \( k \). This is a semiregular bipartite graph with the property that, for any two vertices \( x \) and \( y, \) the number of vertices \( z \) such that \( d \lp x, z \rp = i \) and \( d \lp z, y \rp = j \) depends only on \( d \lp x, y \rp \) and the valency of \( x \). We say that a \textit{distance-biregular} graph is a semiregular bipartite graph with this property. Note that if \( k = \ell, \) this is simply a distance-regular graph.

Distance-regular graphs have a long list of equivalent characterizations, as can be found in such texts as Brouwer, Cohen, and Neumeier~\cite{bcn}. One characterization is that a distance-regular graph is a graph where, for every vertex \( a, \) the distance partition \( \pi_a \) is equitable, and for all vertices \( a, b \), the quotient graphs \( G / \pi_a \) and \( G / \pi_b \) are isomorphic as rooted graphs. Godsil and Shawe-Taylor~\cite{distanceRegularised} studied graphs for which this definition was relaxed slightly.

\begin{thm}[Godsil and Shawe-Taylor~\cite{distanceRegularised}]\label{localDistance}Let \( G \) be a graph with the property that, for every vertex \( a \), the distance partition \( \pi_a \) is equitable. Then \( G \) is either distance-regular or distance-biregular.\end{thm}

In this way, distance-biregular graphs can be seen as a natural relaxation of distance-regular graphs. In addition to~\cite{distanceRegularised}, they have also been studied by Delorme~\cite{delorme}
, though with a different, equivalent, definition. Here, we gather some of the alternate definitions of distance-biregular graphs.

\begin{thm}\label{distanceBiregular}Let \( G = \lp B, C \rp \) be a bipartite graph with diameter \( d \). The following statements are equivalent:
  \begin{enumerate}[label = (\alph*)]
  \item The graph \( G \) is distance-biregular;
  \item For any choice of vertices \( x \) and \( y, \) the number of vertices \( z \) such that \( d \lp x, z \rp = i \) and \( d \lp z, y \rp = j \) depends only on \( d \lp x, y \rp \) and the valency of \( x \);
  \item For every vertex \( a, \) the distance partition \( \pi_a \) is equitable;
  \item Whenever matrix multiplication is well-defined for distance biadjacency matrices, their product is in the span of distance biadjacency matrices.
  \end{enumerate}
\end{thm}

\proof We established $(b)$ as the definition, and from Theorem~\ref{localDistance}, we can see that this is equivalent to $(c)$.

Finally, we will show that $(b)$ and $(d)$ are equivalent. Let \( i, j \) be even integers between 0 and \( d \). Then
\[ \lp N_i N_j \rp_{x,y} = \sum_{z \in B} \lp N_i \rp_{x,z} \lp N_j \rp_{z,y} = \abs{\lsb z \in B: d \lp x, z \rp = i, d \lp y, z \rp = j \rsb}, \]
which is a linear combination of \( N_0, N_1, \ldots, N_d \) precisely when the size of the set is independent of our choice of \( x \) and \( y \). The same argument holds for the other cases in which the multiplication is well-defined.\qed

\section{Orthogonal polynomials}
We can describe the distance biadjacency matrices of the \( \lp k, \ell \rp \)-biregular tree. Let \( B \) be the set of vertices of valency \( k, \) \( C \) the set of vertices of valency \( \ell, \) and \( N \) the biadjacency matrix from \( B \) to \( C \).

Note that \( N_0 = I \) and \( N_1 = N \), and using the block decomposition, we can see that \( NN^T \) counts walks of length 2, so \( N_2 = NN^T - kI. \) Then, for vertices \( b \in B, c \in C \) and some \( i \geq 1, \) we can compute that
\[ \lp N_{2i} N \rp_{b,c} = \sum_{x \in B} \lp N_{2i} \rp_{b,x} N_{x,c} = \sum_{x \sim c} \lp N_{2i} \rp_{b,x}. \]
This entry is equal to the number of neighbours \( c \) has at distance \( 2i \) from \( b, \) which is \( 1 \) if \( d \lp b, c \rp = 2i+1, \) \( \ell -1 \) if \( d \lp b, c \rp = 2i-1, \) and 0 otherwise. Thus, we see that
\begin{equation}\label{nEven}N_{2i} N = N_{2i+1} + \lp \ell - 1 \rp N_{2i-1}.\end{equation}

We can similarly compute that
\begin{equation}\label{nOdd}N_{2i+1} N^T = N_{2i+2} + \lp k-1 \rp N_{2i}.\end{equation}

We can thus define a sequence of polynomials \( \lp F_i^{\lp k, \ell \rp} \lp x, y \rp \rp_{i \geq 0} \) associated to the \( \lp k, \ell \rp \)-biregular tree by making sure they satisfy this same recurrence. If \( N_0, N_1, \ldots, \) are the distance biadjacency matrices of the \( \lp k, \ell \rp \)-biregular tree, then
\[ F_i^{\lp k, \ell \rp} \lp N, N^T \rp = N_i. \]

When evaluating these polynomials using the biadjacency matrix of a finite semiregular bipartite graph, they retain a combinatorial meaning by counting the nonbacktracking walks in the graph.

\begin{lem}\label{backtrackSemiregular}Let \( G = \lp B, C \rp \) be a semiregular graph with valencies \( k \) and \( \ell \), respectively, and biadjacency matrix \( N \) from \( B \) to \( C. \) Then the entries of
  \[ F_{i}^{\lp k, \ell \rp} \lp N, N^T \rp \]
  count the number of nonbacktracking walks of length \( i \) from vertices in \( B \) to vertices in \( B, \) for \( i \) even, and to vertices in \( C \), for \( i \) odd.
\end{lem}

This is analogous to the result in Singleton~\cite{singleton} that the sequence of polynomials associated to the $k$-regular tree count the nonbacktracking walks in a $k$-regular graph. In both cases, the proof follows by induction on the length of the walk \( i\).

In general, when speaking of a semiregular bipartite graph, the associated family of polynomials is understood in terms of the valencies, and so we will suppress the superscripts.

Although this is a function in two variables, we can define a related, single variable, family \( \lp P_i \lp z \rp \rp_{i \geq 0} \) by
\[ P_i \lp z \rp = P_i \lp xy \rp = F_{2i} \lp x, y \rp. \]

Using the recursion for \( N_i \), and by extension \( F_i \lp x, y \rp, \) we see that
\[ P_0 \lp z \rp = 1, \]
\[ P_1 \lp z \rp = z-k, \]
\[ P_2 \lp z \rp = z^2 - \lp 2k+\ell-2 \rp z + k \lp k-1 \rp. \]
For \( i \geq 3, \) we use the recursion in Equations~\ref{nEven} and~\ref{nOdd} to see
\begin{align*}
  F_{2i} \lp \sqrt{z}, \sqrt{z} \rp &= \sqrt{z} F_{2i-1} - \lp k-1 \rp F_{2i-2} \\
  &= \lp z - k +1 \rp F_{2i-2} - \lp \ell-1 \rp \sqrt{z} F_{2i-3} \\
  &= \lp z-k- \ell +2 \rp F_{2i-2} - \lp k-1 \rp \lp \ell-1 \rp F_{2i-4}. \\
\end{align*}
Thus, for \( i \geq 3, \) we have
\begin{equation}\label{pRecursion}P_i \lp z \rp = \lp z -k -\ell+2 \rp P_{i-1} \lp z \rp - \lp k-1 \rp \lp \ell-1 \rp P_{i-2} \lp z \rp.\end{equation}
This satisfies a three-term recurrence, so by Theorem~\ref{favard}, it is a family of monic orthogonal polynomials.

We can define a second family in a similar way. Namely,
\[ I_i \lp xy \rp x = F_{2i+1} \lp x, y \rp. \]

This family satisfies
\[ I_0 \lp z \rp= 1, \]
\[ I_1 \lp z \rp = z- \lp k + \ell -1 \rp. \]
For \( i \geq 2, \) we can use the same expansion as for the \( P_i \) to show that
\begin{equation}\label{iRecursion} I_{i+1} \lp z \rp = z I_i \lp z \rp - \lp \ell + k - 2 \rp I_i \lp z \rp - \lp k-1 \rp \lp \ell -1 \rp I_{i-1} \lp z \rp. \end{equation}
By Theorem~\ref{favard}, it is also a family of monic orthogonal polynomials.

Letting \( N_0, N_1, \ldots \) be the distance biadjacency matrices of the biregular tree, we have that \( P_i \lp N N^T \rp = N_{2i} \) and \( I_i \lp N N^T \rp N = N_{2i+1}. \) Since the biregular tree is distance-biregular, we thus conclude that whenever matrix multiplication is well-defined, the products are in the span of \( P_i \lp N N^T \rp \) or \( I_i \lp N N^T \rp N \), and the coefficients are the intersection numbers of the \( \lp k, \ell \rp \)-biregular tree. By the definition of \( F_i \lp x, y \rp, \) the product of \( P_i \lp xy \rp P_i \lp xy \rp \) will be in the span of \( P_i \lp xy \rp, \) and similarly for the other cases. Letting \( p_s \lp i, j \rp \) denote the intersection numbers of the \( \lp k, \ell \rp \)-biregular tree, we can formalize this discussion to get the following lemma.

\begin{lem}\label{semiregularIntersect}For \( i, j \geq 0, \) there exist positive coefficients \( p_s \lp i, j \rp \) such that
\[ P_i \lp xy \rp P_j \lp xy \rp = \sum_{s=\abs{i-j}}^{i+j} p_{2s} \lp 2i, 2j \rp P_s \lp xy \rp, \]
\[ P_i \lp xy \rp I_j \lp xy \rp x  = \sum_{s=\abs{i-j}}^{i+j} p_{2s+1} \lp 2i, 2j+1 \rp I_s \lp xy \rp x, \]
\[ I_i \lp xy \rp x P_j \lp yx \rp = \sum_{s=\abs{i-j}}^{i+j} p_{2s+1} \lp 2i+1, 2j \rp I_s \lp xy \rp x, \]
and
\[ I_i \lp xy \rp x I_j \lp yx \rp y  = \sum_{s=\abs{i-j}}^{i+j+1} p_{2s} \lp 2i+1, 2j+1 \rp P_s \lp xy \rp. \]
\end{lem}

We can also note that \( p_0 \lp 2i, 2j \rp \) is 0 if \( i \neq j \), and \( p_0 \lp 2i, 2i \rp = 1 \) if \( i =0 \) and \( k \lp \ell-1 \rp^{i} \lp k-1 \rp^{i-1} \) otherwise. Letting \( \delta_{ij} \) be the Kronecker delta, this gives us
\[ p_0 \lp 2i, 2j \rp = P_i \lp k \ell \rp \delta_{i,j}. \]
Similarly,
\[ p_0 \lp 2i+1, 2j+1 \rp = k \lp \ell-1 \rp^{i} \lp k-1 \rp^{i} \delta_{i,j} = k I_i \lp k \ell \rp \delta_{i,j}. \]

\section{Girth and distance-biregular graphs}

We can use the families \( \lp P_i \rp_{i \geq 0} \) and \( \lp I_i \rp_{i \geq 0} \) to obtain another characterization of distance-biregular graphs. From here on, we change our notation. The distance biadjacency matrices from \( B, \) which were previously denoted \( N_0, \ldots, N_d \) are now denoted \( B_0, \ldots, B_d. \) This allows us to also introduce \( C_0, \ldots, C_d \) as the distance biadjacency matrices from \( C \).

\begin{lem}\label{distanceBiregularPolynomial}Let \( G = \lp B, C \rp \) be a bipartite graph with diameter \( d \). Let the distance biadjacency matrices from \( B \) be given by \( B_0, \ldots, B_d, \) and the distance biadjacency matrices from \( C \) be given by \( C_0, \ldots, C_d. \) Then \( G \) is distance-biregular if and only if, for all \( i, \) we have
  \[ P_i \lp B_1 B_1^T \rp, I_i \lp B_1 B_1^T \rp B_1 \in \mathrm{span} \lsb B_0, B_1, \ldots, B_d \rsb \]
  and
  \[ P_i \lp C_1 C_1^T \rp, I_i \lp C_1 C_1^T \rp C_1 \in \mathrm{span} \lsb C_0, C_1, \ldots, C_d \rsb. \]  
\end{lem}

\proof Let \( G \) be a distance-biregular graph. We can use induction to show that \( P_i \lp B_1 B_1^T \rp \) is a linear combination of products of distance biadjacency matrices, so by Theorem~\ref{distanceBiregular}, it is itself in the span of distance biadjacency matrices. The other cases are analogous, which establishes the forward direction.

Now, suppose the other direction. If \( 2i \leq d, \) then by assumption, we can write \( P_i \lp B_1 B_1^T \rp \) as a linear combination of \( B_0, B_2, \ldots, B_{2i} \) and, since by Lemma~\ref{backtrackSemiregular} we know \( P_i \lp B_1 B_1^T \rp \) counts nonbacktracking walks of length \( 2i \), the coefficient of \( B_{2i} \) is necessarily nonzero. It follows by induction that we can write \( B_{2i} \) as a polynomial in \( B_1 B_1^T \), and therefore we can write it as a linear combination of polynomials in \( \lp P_i \lp B_1 B_1^T \rp \rp_{i \geq 0} \).

\vspace{.1in}

We can similarly show that \( B_{2i+1}, C_{2i}, \) and \( C_{2i+1} \) are linear combinations of \( \lp I_i \lp B_1 B_1^T \rp B_1 \rp_{i \geq 0}, \) \( \lp P_i \lp C_1 C_1^T \rp \rp_{i \geq 0} \), and \( \lp I_i \lp C_1 C_1^T \rp C_1 \rp_{i \geq 0} \), so applying Lemma~\ref{semiregularIntersect} and Theorem~\ref{distanceBiregular}, we conclude that \( G \) is distance-biregular.\qed

Let \( G \) be a regular graph with girth \( g \), diameter \( d, \) and \( d+1 \) distinct eigenvalues. Haemers and van Dam~\cite{vanDam} showed that if \( g \geq 2d-1 \), then \( G \) is distance-regular, and Abiad, van Dam, and Fiol~\cite{abiad} showed that for a bipartite graph \( G \), if \( g \geq 2d-2, \) then \( G \) is distance-regular. A new result shows that an analogous statement is true for distance-biregular graphs.

\begin{thm}\label{girthBiregular}Let \( G = \lp B, C \rp \) be a connected \( \lp k, \ell \rp \)-semiregular bipartite graph with girth \( g \), diameter \( d, \) and \( d+1 \) distinct eigenvalues. If \( g \geq 2d-2, \) then \( G \) is distance-biregular.
\end{thm}

\proof Let \( B_0, B_1, \ldots, B_{d} \) be the distance biadjacency matrices from \( B \). Let \( e \) be the largest integer such that \( 2e \leq d.  \) We wish to show that, for all \( i, \)
\[ P_i \lp B_1 B_1^T \rp \in \mathrm{span} \lsb B_0, B_2, \ldots, B_{2e} \rsb. \]

By Lemma~\ref{backtrackSemiregular}, \( P_i \lp B_1 B_1^T \rp \) counts the nonbacktracking walks of \( G, \) and since \( g \geq 2d-2, \) we know that for all \( i \leq d-2, \) any nonbacktracking walk in \( G \) must be a path, and there is at most one such path. From this, we conclude that, for all \( i \in \lsb 0, 1, \ldots, e-1 \rsb, \) we have
\[ P_i \lp B_1 B_1^T \rp = B_{2i}. \]

Let \( J_{B,B} \) be the \( \abs{B} \times \abs{B} \) all-ones matrix, and define \( J_{B,C}, J_{C,B} \) and \( J_{C,C} \) similarly. Let \( E_1 \) be the spectral idempotent corresponding to the eigenvalue \( \sqrt{k \ell} \), and note
\[ \prod_{r=2}^d \frac{\lp x- \theta_r \rp}{\sqrt{k \ell}-\theta_r} = E_1 = \pmat{\frac{1}{2 \abs{B}} J_{B,B} & \frac{\sqrt{\ell}}{2 \abs{B} \sqrt{k}} J_{B,C} \\ \frac{\sqrt{k}}{2 \abs{C} \sqrt{\ell}} J_{C,B} & \frac{1}{2 \abs{C}} J_{C,C}}. \]
This is a polynomial of degree \( d \) which, evaluated at \( A \), gives \( E_1. \) It follows that there must be a polynomial of degree \( e \) which, evaluated at \( NN^T, \) gives \( J_{B,B}. \) Thus, there exist scalars \( a_i \) with \( a_e \neq 0 \) such that
\[ J_{B,B} = \sum_{i=0}^e a_i P_i \lp B_1 B_1^T \rp = \sum_{i=0}^e B_{2i}, \]
from which we get
\[ P_e \lp B_1 B_1^T \rp = \frac{1}{a_e} \lp B_{2e} + \sum_{i=0}^{e-1} \lp 1-a_i \rp B_{2i} \rp \in \mathrm{span} \lsb B_0, B_2, \ldots, B_{2e} \rsb. \]

For \( i \geq e, \) we can proceed by induction to show that \( P_{i+1} \lp B_1 B_1^T \rp \in \mathrm{span} \lsb B_0, B_2, \ldots, B_{2e} \rsb. \) Using Equation~\ref{pRecursion}, it suffices to show
\[ B_1B_1^T P_i \lp B_1 B_1^T \rp \in \mathrm{span} \lsb B_0, B_2, \ldots, B_{2e} \rsb. \]
From the inductive hypothesis, we can write
\[ P_i \lp B_1 B_1^T \rp = \sum_{i=0}^e b_i B_{2i} =  b_e B_{2e} + \sum_{i=0}^{e-1} b_i P_i \lp B_1 B_1^T \rp. \]

Note that
\[ B_1B_1^T B_e = B_1B_1^T J_{B,B} - B_1B_1^T \sum_{i=0}^{e-1} P_i \lp B_1 B_1^T \rp, \]
and using Equation~\ref{pRecursion} again, we have shown that
\[ B_1B_1^T \sum_{i=0}^{e-1} P_i \lp B_1 B_1^T \rp \in \mathrm{span} \lsb P_0, P_1, \ldots, P_e \rsb. \]
Thus, we have that \( P_{i+1} \lp B_1 B_1^T \rp \in \mathrm{span} \lsb B_0, B_2, \ldots, B_{2e} \rsb \).

Since the recursion for \( P_i \) and \( I_i \) are the same, it suffices to prove the base case for \( I_i \lp B_1 B_1^T \rp B_1. \) Let \( e' \) be the largest integer such that \( 2e'+1 \leq d. \) Summing up all the odd distance biadjacency matrices, we get the \( \abs{B} \times \abs{C} \) matrix of all ones. Since \( J_{B,B} \) is polynomial in \( B_1 B_1^T \), we may write
\[ \ell \sum_{i=0}^{e'} B_{2i+1}  = J_{B,B} B_1 = \sum_{i=0}^{e'} b_i I_i \lp B_1 B_1^T \rp B_1, \]
for some scalars \( b_0, \ldots, b_{e'} \) with \( b_{e'} \neq 0. \) This gives us
\[ I_{e'} \lp B_1 B_1^T \rp B_1 = \frac{1}{b_{e'}} \lp \ell B_{2e'+1} + \sum_{i=0}^{e'-1} \lp \ell-b_i \rp B_{2i+1} \rp. \]

By induction, for all \( i \geq 0, \) we have \( I_i \lp B_1 B_1^T \rp B_1 \in \mathrm{span} \lsb B_1, B_3, \ldots, B_{2e'+1} \rsb. \) The same arguments apply to \( P_i \lp B_1^T B_1 \rp \) and \( I_i \lp B_1^T B_1 \rp B_1^T, \) so by Lemma~\ref{distanceBiregularPolynomial}, we conclude that \( G \) is distance-biregular.\qed

\section{Linear programming bound for bipartite semiregular graphs}

It is a well-known fact that the largest eigenvalue is at most the maximum valency of a graph, with equality holding if and only if the graph is regular. Similarly, the largest eigenvalue of a bipartite graph is at most \( \sqrt{\Delta_1 \Delta_2}, \) where \( \Delta_1 \) is the maximum valency of vertices in the first partition, and \( \Delta_2 \) is the maximum valency of vertices in the second partition. Equality holds if and only if the graph is semiregular. These results can be found, for instance, on page 173 of Godsil and Royle~\cite{yellow}.

Since the largest eigenvalue of a bipartite semiregular graph is determined by the valencies, we turn our attention to the second-largest eigenvalue. In particular, we generalize the linear programming bound in Cioab\u{a}, Koolen, and Nozaki~\cite{bipartite} to the semiregular case.

\begin{thm}\label{lpSemiregular}Let \( G = \lp B, C \rp \) be a semiregular graph with valencies \(k, \ell. \) Let \( \sqrt{k \ell} = \theta_1 > \theta_2 > \cdots > \theta_d \) be the set of distinct, nonnegative eigenvalues of \( G \). Suppose there exists a polynomial
  \[ f \lp z \rp = \sum_{i=0}^t f_i P_i \lp z \rp \]
  such that
  \begin{itemize}
  \item For eigenvalues \( \theta_r \neq \theta_1, \) we have \( f \lp \theta_r^2 \rp \leq 0 \);
  \item The coefficient \( f_0 \) is positive;
  \item For \( 1 \leq i \leq t \), the coefficient \( f_i \) is nonnegative. 
  \end{itemize}
  Then
  \[ \abs{B} \leq \frac{f \lp k \ell \rp}{f_0}. \]
  If equality holds, then \( f \lp \theta_r^2 \rp = 0 \) for all \( \theta_r \neq \pm \sqrt{k \ell}. \) Further, if equality holds and \( f_i > 0 \) for each \( i \in \lsb 1, \ldots, t \rsb \), then \( G \) has girth at least \( 2t+2. \)
\end{thm}      

\proof Let \( A \) be the adjacency matrix of \( G \), and \( N \) be the biadjacency matrix from \( B \) to \( C \). Since \( NN^T \) and \( N^T N \) share nonzero eigenvalues, we can take the block decomposition of \( A \) to show that every eigenvalue of \( N N^T \) is an eigenvalue of \( A^2 \), and every nonzero eigenvalue of \( A^2 \) is an eigenvalue of \( N N^T \). Thus, we can write the spectral decomposition
\[ NN^T = \sum_{r=1}^{d'} \theta_r^2 E_r, \]
where \( d' \) is equal to either \( d \) or \( d-1, \) \( \theta_r^2 \) is an eigenvalue of \( NN^T, \) and \( E_r \) are the spectral idempotents for \( NN^T \).

Then we have
\[ f \lp k \ell \rp E_1 + \sum_{r=2}^d f \lp \theta_r^2 \rp E_r = f \lp N N^T \rp = f_0 I + \sum_{i=1}^t f_i P_i \lp N N^T \rp. \]

Taking the trace of both sides gives us
\[  f \lp k \ell \rp \geq f \lp k \ell \rp + \sum_{r=2}^d f \lp \theta_r^2 \rp \mathrm{tr} \lp E_r \rp = f_0 \abs{B} + \sum_{i=1}^t f_i \mathrm{tr} \lp P_i \lp NN^T \rp \rp \geq f_0 \abs{B}. \]

When equality holds, we must have that for all \( r \neq 1, \) with the possible exception of \( \theta_d = 0, \) we have \( f \lp \theta_r^2 \rp = 0 \). Further, for all \( i \geq 0, \) we have
\[ f_i \mathrm{tr} \lp P_i \lp NN^T \rp \rp = 0, \]
and, by Lemma~\ref{backtrackSemiregular}, if \( f_i > 0 \) for all \( i \in \lsb 1, \ldots, t \rsb, \) we must have that the girth is at least \( 2t +2. \)\qed

\section{Upper bound for semiregular bipartite graphs with given second eigenvalue}

Let \( b \lp k, \ell, \theta \rp \) be the maximum number of vertices of valency \( k \) in a \( \lp k, \ell \rp \)-semiregular graph with second eigenvalue at most \( \theta \). We wish to use our bound from the previous section to find an upper bound on \( b \lp k, \ell, \theta \rp. \)

We will need another family of orthogonal polynomials, called \( \lp J \lp z \rp \rp_{i \geq 0}, \) defined by
\begin{equation}\label{ji} J_i \lp z \rp = \sum_{j=0}^{i} I_j \lp z \rp.\end{equation}
It can be verified by induction that, for \( i \geq 2, \) the family \( J_i \lp z \rp \) satisfies the same recursion as \( P_i \lp z \rp \) and \( I_i \lp z \rp, \) so by Theorem~\ref{favard}, it forms a family of monic orthogonal polynomials. Note that
\[ I_i = J_i - J_{i-1}. \]

Let \( B \lp k, \ell, 2t+1, c \rp \) be the \( \lp 2t+1 \rp \times \lp 2t+1 \rp \) tridiagonal matrix with lower diagonal \( \lp 1, \ldots, 1, c, \ell \rp \), zero along the main diagonal, and row sum alternating between \( k \) and \( \ell. \)

\begin{thm}\label{odd}If \( \theta \) is the second-largest eigenvalue of \( B \lp k, \ell, 2t+1, c \rp \), then
\[ b \lp k, \ell, \theta \rp \leq \ell \sum_{i=0}^{t-2} \lp \ell-1 \rp^i \lp k - 1 \rp^{i} + \frac{\ell \lp \ell -1 \rp^{t-1} \lp k-1 \rp^{t-1}}{c}. \]
\end{thm}

\proof Consider the matrix \( xI - B \lp k, \ell, 2t+1, c \rp \). 

We can verify by induction that, for \(i \leq 2t-1, \) the determinant of the principal \( i \times i \) matrix formed by the first \( i \) rows and \( i \) columns of this matrix is \( F_i \lp x, x \rp. \) Thus expanding along the bottom two rows and right column of \( xI - B \lp k, \ell, 2t+1, c \rp \) gives us that the determinant is
\[ x^2 F_{2t-1} \lp x, x \rp - c \lp k-1 \rp x F_{2t-2} \lp x, x \rp - k \lp \ell-c \rp F_{2t-1} \lp x, x \rp. \]
Using Equation~\ref{nOdd} and the definitions of \( F_i \lp x, x \rp \) and \( I \lp x \rp, \) this is equivalent to
\[ x \lp x^2 I_{t-1} \lp x^2 \rp - c \lp k-1 \rp \lp I_{t-1} \lp x^2 \rp - \lp \ell-1 \rp I_{t-2} \lp x^2 \rp \rp - k \lp \ell-c \rp I_{t-1} \lp x \rp \rp, \]
which simplifies to
\begin{equation}\label{det1} x \lp \lp x^2 - k \ell \rp I_{t-1} \lp x^2 \rp + c \lp I_{t-1} \lp x^2 \rp - \lp k-1 \rp \lp \ell-1 \rp I_{t-2} \lp x^2 \rp \rp \rp.\end{equation}

Using Equation~\ref{iRecursion}, we compute that
\[ I_{t-1} \lp x^2 \rp - \lp k-1 \rp \lp \ell-1 \rp I_{t-2} \lp x^2 \rp = \lp x^2 - k \ell + 1 \rp I_{t-2} \lp x^2 \rp - \lp k-1 \rp \lp \ell-1 \rp I_{t-3} \lp x^2 \rp, \]
so we can inductively rewrite Equation~\ref{det1} as
\[ x \lp x^2 - k \ell \rp \lp I_{t-1} \lp x^2 \rp + c \sum_{i=0}^{t-2} I_i \lp x^2 \rp \rp. \]
Using Equation~\ref{ji}, this simplifies to
\[ x \lp x^2 - k \ell \rp \lp \lp c -1 \rp J_{t-2} \lp x^2 \rp + J_{t-1} \lp x^2 \rp \rp. \]

The eigenvalues of \( B \lp k, \ell, 2t+1, c \rp \) are the roots of \( \det \lp xI - B \lp k, \ell, 2t+1, c \rp \rp \). Therefore, the nontrivial, nonzero eigenvalues of \( B \lp k, \ell, 2t+1, c \rp \) are the square roots of the roots of the expression
\[ \lp c -1 \rp J_{t-2} \lp z \rp + J_{t-1} \lp z \rp. \]
Let \( \sqrt{k \ell} =  \theta_1 > \theta_2 > \cdots > \theta_t \) be the positive eigenvalues of \( B \lp k, \ell, 2t+1, c \rp. \) Then
\[ \prod_{r=2}^{t} \lp x- \theta_r^2 \rp = \lp c -1 \rp J_{t-2} \lp z \rp + J_{t-1} \lp z \rp. \]

Therefore, we can define
\begin{align*}
  f \lp z \rp &= z \lp z- \theta_2^2 \rp \prod_{r=3}^t \lp z-\theta_r^2 \rp^2 \\
  &= \frac{\lp c -1 \rp J_{t-2} \lp z \rp + J_{t-1} \lp z \rp}{z-\theta^2} \lp cz \sum_{i=0}^{t-2} I_{i} \lp z \rp + z I_{t-1} \lp z \rp \rp.
\end{align*}

For all \( r \in \lsb 1, \ldots, d \rsb, \) we have \( f \lp \theta_r^2 \rp \leq 0. \) Since \( f \lp z \rp \) is a polynomial in \( z \) of degree \( 2t-2, \) we may write it as
\[ f \lp z \rp = \sum_{i=0}^{2t-2} f_i P_i \lp z \rp, \]
and we will show that \( f_i > 0 \) to apply the linear programming bound.

By Proposition~\ref{posCoeff} and the definition of \( J_i \), we know there exists positive coefficients \( j_0, \ldots, j_{t-2} \) such that
\begin{equation}\label{ji} \frac{\lp c -1 \rp J_{t-2} \lp z \rp + J_{t-1} \lp z \rp}{z-\theta^2} = \sum_{i=0}^{t-2} j_i I_i \lp z \rp.\end{equation}
Thus we may write
\[ f \lp z \rp =  \lp \sum_{i=0}^{t-2} j_i I_i \lp z \rp \rp \lp cz \sum_{i=0}^{t-2} I_{i} \lp z \rp + z I_{t-1} \lp z \rp \rp \]
for \( j_i > 0. \)

By Lemma~\ref{semiregularIntersect}, we get that \( f_i > 0 \) for all \( i=0, \ldots, 2t-2. \) Further, from the discussion after Lemma~\ref{semiregularIntersect} and Equation~\ref{ji}, we have
\[ f_0 = ck \sum_{i=0}^{t-2} j_i I_i \lp k \ell \rp = ck \frac{\lp c -1 \rp J_{t-2} \lp k \ell \rp + J_{t-1} \lp k \ell \rp}{k \ell-\theta^2} \]

Applying Theorem~\ref{lpSemiregular}, we then get
\begin{align*}
  \abs{B} &\leq \frac{f \lp k \ell \rp}{f_0} \\
  &= \ell \sum_{i=0}^{t-2} I_{i} \lp k \ell \rp + \ell \frac{I_{t-1} \lp k \ell \rp}{c} \\
  &= \ell \sum_{i=0}^{t-2} \lp \ell-1 \rp^i \lp k - 1 \rp^{i} + \frac{\ell \lp \ell -1 \rp^{t-1} \lp k-1 \rp^{t-1}}{c}.\tag*{\sqr53}
\end{align*}

A close variation of this proof gives us an analogue for even matrices. Let \( B \lp k, \ell, 2t, c \rp \) be the \( 2t \times 2t \) tridiagonal matrix with lower diagonal \( \lp 1, \ldots, 1, c, k \rp \), zero along the main diagonal, and row sum alternating between \( k \) and \( \ell. \) For this proof, we also need a new family of orthogonal polynomials, called \( \lp Q \lp z \rp \rp_{i \geq 0} \), defined by
\[ Q_i \lp z \rp = \sum_{j=0}^{i} P_i \lp z \rp. \]

\begin{thm}\label{even}If \( \theta \) is the second-largest eigenvalue of \( B \lp k, \ell, 2t, c \rp \), then
  \[ b \lp k, \ell, \theta \rp \leq 1 +  k \sum_{i=1}^{t-2} \lp \ell-1 \rp^i \lp k - 1 \rp^{i-1} + \frac{k \lp \ell -1 \rp^{t-1} \lp k-1 \rp^{t-2}}{c}. \]
\end{thm}

\proof Apply Theorem~\ref{lpSemiregular} with
\[ f \lp z \rp = \frac{\lp \lp c -1 \rp Q_{t-2} \lp z \rp + Q_{t-1} \lp z \rp \rp^2}{z-\theta^2}.\qed \]

It is worth noting that although Theorem~\ref{even} has the same form as the bounds found by Cioab\u{a}, Koolen, Nozaki, and Vermette~\cite{cioaba1}, and Cioab\u{a}, Koolen, Mimura, Nozaki, and Okuda~\cite{hypergraph}, the choice of \( c \) might vary depending on whether the bound comes from considering a semiregular bipartite graph, a uniform and regular hypergraph, or the regular halved graph of a semiregular bipartite graph.

\begin{lem}\label{biregularEquality}Equality holds in Theorems~\ref{odd} and~\ref{even} if and only if there exists a distance-biregular graph whose quotient matrix with respect to the distance-partition from a vertex of valency \( k \) is \( B \lp k, \ell, 2t+1, c \rp \) or \( B \lp k, \ell, 2t, c \rp \), respectively.\end{lem}

\proof If \( G \) is a distance-biregular graph with quotient matrix of the form \( B \lp k, \ell, 2t+1, c \rp \) or \( B \lp k, \ell, 2t, c \rp \), then a simple counting argument reveals \( G \) has this maximal number of vertices.

If equality holds in Theorem~\ref{odd}, then it must hold in Theorem~\ref{lpSemiregular}.

From the proof of Theorem~\ref{odd}, we recall that
\[ f \lp z \rp = \sum_{i=0}^{2t-2} f_i P_i \lp z \rp, \]
where \( f_i > 0 \), so when we apply Theorem~\ref{lpSemiregular}, we see that \( G \) has girth at least \( 4t-2 \). This implies that the diameter \( d \geq 2t-1. \) If \( d = 2t-1, \) then it must be a generalized polygon, which is distance-biregular, as discussed by Godsil and Shawe-Taylor~\cite{distanceRegularised}. Otherwise, we have that \( d \geq 2t. \)

From Theorem~\ref{lpSemiregular}, we also see that every nonzero, nontrivial eigenvalue of \( G \) is a square root of a zero of a \( t-1 \) degree polynomial. It follows that \( G \) must have at most \( 2t+1 \) distinct eigenvalues, and thus \( G \) has diameter \( d = 2t. \) 

Therefore
\[ g \geq 4t -2 = 2d-2, \]
so by Theorem~\ref{girthBiregular}, \( G \) must be distance-biregular.

Similarly, a graph obtaining the bound in Theorem~\ref{even} has girth \( g \geq 4t-4 \) and at most \( 2t-2 \) distinct eigenvalues, so either \( G \) is a generalized polygon, or \( G \) has diameter \( d=2t-1 \) and girth
\[ g \geq 4t-4 = 2d -2. \] 
In either case, \( G \) must be distance-biregular.\qed

We would like to know that, when \( b \lp k, \ell, \theta \rp \) is well-defined, we can in fact find a choice of matrix \( B \lp k, \ell, t, c \rp \) to apply either Theorem~\ref{odd} or Theorem~\ref{even} to find an upper bound. Feng and Li~\cite{fengLi} proved a generalization of the Alon-Boppana theorem for \( \lp k, \ell \rp \)-semiregular bipartite graphs. If \( 0 < \theta < \sqrt{k-1} + \sqrt{\ell-1}, \) then there are only finitely many graphs that have second-largest eigenvalue at most \( \theta \). Thus, it is precisely this range of eigenvalues for which we would like to find \( B \lp k, \ell, \theta \rp. \)

\begin{prop}Let \( k, \ell \geq 3 \) and let \( 0 < \theta < \sqrt{k-1}+\sqrt{\ell-1} \). Then there exists a matrix \( B \lp k, \ell, t, c \rp \) with second-largest eigenvalue \( \theta \).\end{prop}

\proof Let \( T \) be the quotient matrix of the \( \lp k, \ell \rp \)-biregular tree, and, for \( i \geq 1, \) let \( T_i \) be the principal submatrix formed by the first \( i \) rows and columns. As noted in the proof of Theorem~\ref{odd}, we can show by induction that
\[ \mathrm{det} \lp xI - T_i \rp = F_i \lp x, x \rp. \]
Thus, the eigenvalues of \( T_i \) are the roots of \( F_i \lp x, x \rp. \)

On the other hand, as \( i \) tends towards infinity, the sequence \( \lp T_i \rp_{i \geq 0} \) converges to \( T, \) and the quotient matrix of a distance-biregular graph shares the same eigenvalues as the graph. From Godsil and Mohar~\cite{infiniteSpectrum}, we know that the \( \lp k, \ell \rp \)-biregular tree has second-largest eigenvalue \( \sqrt{k-1}+\sqrt{\ell-1} \). Therefore, the second-largest zero of \( F_i \lp x, x \rp \) converges to \( \sqrt{k-1} + \sqrt{\ell-1}. \)

Thus, there exists some \( t \) such that \( P_t \lp \theta^2 \rp < 0, \) but \( P_{t-1} \lp \theta^2 \rp \geq 0 \) or \( I_t \lp \theta^2 \rp < 0 \) but \( I_{t-1} \lp \theta^2 \rp \geq 0. \) The two cases are analogous, so we will prove the second case.

Let
\[ c= \frac{-I_t \lp \theta^2 \rp}{J_{t-1} \lp \theta^2 \rp}. \]
Note that by construction, \( c > 0 \).

We have that
\begin{align*}
  J_{t} \lp \theta^2 \rp + \lp c-1 \rp J_{t-1} \lp \theta^2\rp &= J_t \lp \theta^2 \rp - J_{t-1} \lp \theta^2 \rp + c J_{t-1} \lp \theta^2 \rp \\
  &= I_t \lp \theta^2 \rp - \frac{I_{t} \lp \theta^2 \rp}{J_{t-1} \lp \theta^2 \rp} J_{t-1} \lp \theta^2 \rp \\
  &= 0,
\end{align*}
so \( \theta \) is a root of \( J_t \lp z^2 \rp + \lp c-1 \rp J_{t-1} \lp z^2 \rp \), and thus an eigenvalue of \( B \lp k, \ell, 2t+3, c \rp. \)

It remains to show that \( \theta \) is in fact the second-largest eigenvalue. From Equation~\ref{iRecursion}, we compute
\begin{align*}
  I_2 \lp x \rp - \lp k-1 \rp \lp \ell-1 \rp I_1 \lp x \rp &= \lp x - k \ell \rp I_1 \lp x \rp + x - k - \ell + 1 - \lp k-1 \rp \lp \ell-1 \rp \\
  &= \lp x-k \ell \rp I_1 + x - k \ell \\
  &= \lp x-k \rp J_1 \lp x \rp.
\end{align*}
Then we may inductively assume that, for \( i \geq 1, \) we have
\[ I_i \lp x \rp - \lp k-1 \rp \lp \ell-1 \rp I_{i-1} = \lp x- k \ell \rp J_{i-1} \lp x \rp. \]
Expanding out using Equation~\ref{iRecursion}, we show that
\begin{align*}
  I_{i+1} \lp x \rp - \lp k-1 \rp \lp \ell-1 \rp I_i \lp x \rp &= \lp x - k \ell \rp I_i \lp x \rp + I_i \lp x \rp - \lp k-1 \rp \lp \ell-1 \rp I_{i-1} \lp x \rp \\
  &= \lp x-k \ell \rp I_i \lp x \rp + \lp x-k \ell \rp J_{i-1} \lp x \rp \\
  &= \lp x-k \ell \rp J_i \lp x \rp.
\end{align*}
Therefore, by Lemma~\ref{interlace}, the nontrivial eigenvalues of \( I_{t-1} \lp x^2 \rp \) interlace the zeroes of \( J_{t-1} \lp x^2 \rp \). In turn, the zeros of \( J_{t-1} \lp x^2 \rp \) interlace the zeros of
\[ J_{t-1} \lp x^2 \rp + \lp c-1 \rp J_{t} \lp x^2 \rp. \]

In particular, if \( \theta \) is not the second-largest eigenvalue of \( B \lp k, \ell, 2t+2, c \rp \), then it must be strictly less than the second-largest zero of \( I_{t-1} \lp x^2 \rp, \) which contradicts our choice of \( t \). Thus we have found a matrix \( B \lp k, \ell, 2t+3, c \rp \) with second-largest eigenvalue \( \theta. \)\qed

\section{Comparison to previous bounds}
A \textit{partial geometry} pg \( \lp \ell, k, \alpha \rp \) is an incidence structure in which there is at most one line between any two points, each line is incident to \( \ell \) points, each point is incident to \( k \) lines, and if \( p \) and \( b \) are not incident, then there are \( \alpha \) coincident pairs \( \lp q, c \rp \) such that \( q \) is incident to \( b \) and \( c \) is incident to \( p \).

Infinite families of partial geometries exist, such as the construction by de Clerck, Dye, and Thas~\cite{clerck} of partial geometries pg \( \lp 2^{2n-1}, 2^{2n-1}+1, 2^{2n-2} \rp \) for \( n \geq 1. \) Mathon~\cite{mathon} constructed an infinite family of partial geometries pg \( \lp 3^{2n}, \frac{1}{2} \lp 3^{4n}+1 \rp, \frac{1}{2} \lp 3^{2n}-1 \rp \rp \) for \( n \geq 1.\)

\begin{ex}Let \( 1 < \alpha < \ell \) be an integer such that a partial geometry \(\text{pg }\lp \ell, k, \alpha \rp \) exists with \( k, \ell \geq 2. \) Viewing this as a semiregular bipartite graph, and starting from a point, we get the distance partition
  \[ \pmat{0 & k \\ 1 & 0 & \ell-1 \\ & 1 & 0 & k -1 \\ & & \alpha & 0 & \ell-\alpha \\ & & & k & 0}, \]
  and second eigenvalue \( \sqrt{k+\ell-\alpha-1}. \)

  Applying the bound of Theorem~\ref{odd}, we get that
  \[ b \lp k, \ell, \sqrt{k+\ell-\alpha-1} \rp \leq \ell + \frac{\ell \lp \ell-1 \rp \lp k-1 \rp}{\alpha} = \frac{\ell \lp \lp \ell-1 \rp \lp k-1 \rp + \alpha \rp}{\alpha}, \]
  which is precisely the number of points in a partial geometry, so this bound is tight.
\end{ex}

By Lemma~\ref{biregularEquality}, we know that the bounds are tight when dealing with a distance-biregular graph. The bound gives us the number of vertices in a partition, and we know that the halved graphs of a distance-biregular graph are distance-regular. Thus it is natural to ask whether Theorem~\ref{odd} and Theorem~\ref{even} give us an improvement over merely applying the bounds for regular graphs found in Cioab\u{a}, Koolen, Nozaki, and Vermette~\cite{cioaba1} to the halved graphs.

\begin{ex}Let \( 1 < \alpha < \ell \) be an integer such that a partial geometry \(\text{pg }\lp \ell, k, \alpha \rp \) exists with \( k, \ell \geq 2. \) The point graph has distance partition
  \[ \pmat{0 & k \lp \ell-1 \rp \\ 1 & \ell - k - 1 + \alpha \lp k-1 \rp & \lp k-1 \rp \lp \ell - \alpha \rp \\ & k \alpha & k \lp \ell-1 -\alpha \rp}, \]
  and second eigenvalue \( \ell-1-\alpha. \) Using the bound in~\cite{cioaba1}, we see the number of points is at most
  \[ b \lp k \lp \ell-1 \rp, \ell-1-\alpha \rp \leq 1 + k \lp \ell-1 \rp + \frac{k \lp \ell-1 \rp \lp k \lp \ell-1 \rp -1 \rp \lp \ell- \alpha \rp}{k \lp \ell-1 \rp -1 - \lp \ell - 1 - \alpha \rp^2}. \]
  Note that, since \( 2 \leq \alpha \leq \ell-1, \) and \( 2 \leq k, \) we have
  \[ k \lp \alpha -1 \rp \lp k \lp \ell-1 \rp -1 \rp + \lp k \lp \ell-1 \rp -1 \rp+ \lp k-1 \rp \lp \ell-1 - \alpha \rp^2 > 0. \]
  Therefore,
  \[ \alpha k \lp k \lp \ell-1 \rp -1 \rp > \lp k-1 \rp \lp k \lp \ell-1 \rp -1  - \lp \ell-1-\alpha \rp^2 \rp, \]
  which we can rewrite as
  \[ \frac{k \lp k \lp \ell-1 \rp -1 \rp}{k \lp \ell-1 \rp -1 - \lp \ell-1 - \alpha \rp^2} > \frac{k-1}{\alpha}. \]
  Multiplying both sides by \( \lp \ell-1 \rp \lp \ell-\alpha \rp \) gives us
  \[ \frac{k \lp \ell-1 \rp \lp k \lp \ell-1 \rp -1 \rp \lp \ell-\alpha \rp}{k \lp \ell-1 \rp -1 - \lp \ell-1 - \alpha \rp^2} > \frac{\lp k-1 \rp \lp \ell-1 \rp \lp \ell-\alpha \rp}{\alpha}. \]

  We compute that
  \[ 1 + k \lp \ell-1 \rp + \frac{\lp k-1 \rp \lp \ell-1 \rp \lp \ell-\alpha \rp}{\alpha} = \frac{\ell \lp \lp k-1 \rp \lp \ell -1 \rp +\alpha \rp}{\alpha}.\]
  This in turn implies that
  \[ 1 + k \lp \ell-1 \rp + \frac{k \lp \ell-1 \rp \lp k \lp \ell-1 \rp -1 \rp \lp \ell- \alpha \rp}{k \lp \ell-1 \rp -1 - \lp \ell - 1 - \alpha \rp^2} > \frac{\ell \lp \ell-1 \rp \lp k-1 \rp + \alpha}{\alpha}, \]
  and therefore, the bound given in~\cite{cioaba1} by considering the point graph is not tight.
\end{ex}

In another paper, Cioab\u{a}, Koolen, Mimura, Nozaki, and Okuda~\cite{hypergraph} proved a bound for hypergraphs, with the comment that the bound also applies to semiregular bipartite graphs. The same example of a partial geometry can show graphs where those hypergraph bounds are not tight, even though the semiregular bipartite bound described here is.

\begin{ex}Let \( 1 < \alpha < \ell \) be an integer such that a partial geometry \(\text{pg }\lp \ell, k, \alpha \rp \) exists with \( k, \ell \geq 2. \) The point graph is the same as before, with second eigenvalue \( \ell-1-\alpha. \) Using the bound from~\cite{hypergraph}, we have that the number of points is bounded above by
  \[ 1 + k \lp \ell-1 \rp + \frac{k \lp k-1 \rp \lp \ell-1 \rp^2 \lp \ell- \alpha \rp}{k \lp \ell-1 \rp + \lp \alpha-1 \rp \lp \ell-1-\alpha \rp}. \]
  
  Since \( \alpha, k, \ell \geq 2, \) we have that
  \[ \lp \alpha-1 \rp \lp k-1 \rp \lp \ell-1 \rp + \alpha \lp \alpha-1 \rp > 0. \]
  Therefore,
  \[ \alpha k \lp \ell-1 \rp > k \lp \ell-1 \rp + \lp \alpha-1 \rp \lp \ell-1 - \alpha \rp, \]
  which we can rewrite as
  \[ \frac{k \lp \ell-1 \rp}{k \lp \ell-1 \rp + \lp \alpha-1 \rp \lp \ell-1 - \alpha \rp} > \frac{1}{\alpha}. \]
  Multiplying both sides by \( \lp k-1 \rp \lp \ell-1 \rp \lp \ell-\alpha \rp \) gives us
  \[ \frac{k \lp k-1 \rp \lp \ell-1 \rp^2 \lp \ell-\alpha \rp}{k \lp \ell-1 \rp + \lp \alpha-1 \rp \lp \ell-1 - \alpha \rp} > \frac{\lp k-1 \rp \lp \ell-1 \rp \lp \ell-\alpha \rp}{\alpha}. \]

  We compute that
  \[ 1 + k \lp \ell-1 \rp + \frac{\lp k-1 \rp \lp \ell-1 \rp \lp \ell-\alpha \rp}{\alpha} = \frac{\ell \lp \lp k -1 \rp \lp \ell-1 \rp + \alpha \rp}{\alpha}. \]
  It follows that
  \[  1 + k \lp \ell-1 \rp + \frac{k \lp k-1 \rp \lp \ell-1 \rp^2 \lp \ell- \alpha \rp}{k \lp \ell-1 \rp + \lp \alpha-1 \rp \lp \ell-1-\alpha \rp} \]
  and therefore, the bound given in~\cite{hypergraph} is also not tight for partial geometries.
\end{ex}

Thus, using constructions of infinite families such as those by de Clerck, Dye, and Thas~\cite{clerck} or Mathon~\cite{mathon}, we have infinite families of distance-biregular graphs where the bound in Theorem~\ref{odd} is tight, but the bounds by Cioab\u{a}, Koolen, Nozaki, and Vermette in~\cite{cioaba1} and Cioab\u{a}, Koolen, Mimura, Nozaki, and Okuda in~\cite{hypergraph} are not.

\section*{Acknowledgements}
I would like to thank Chris Godsil for his feedback and suggestions at every stage of the process. I would also like to thank Sebastian Cioab\u{a} and anonymous referees for their helpful comments on earlier drafts of the paper.


\end{document}